\documentclass{article}
\usepackage{amsmath,amsthm,amssymb,latexsym, amsfonts,}

\theoremstyle{plain}
\newtheorem{teor}{Theorem}[section]
\newtheorem{defin}[teor]{Definition}
\newtheorem{propo}[teor]{Proposition} 
\newtheorem{obs}[teor]{Remark} 
\newtheorem{theo}[teor]{Theorem}
 
\newtheorem{example}[teor]{Example} 
\newtheorem{coro}[teor]{Corollary}

\def\proof{{\noindent \bf Proof:} \hspace{0.1 cm}}
\newcommand{\cqd}{\hfill \rule{2mm}{2mm}\vspace{.3cm}}
\newcommand{\dis}{\displaystyle}

\def \Z {\mathbb Z}

\def \N {\mathbb N}

\def \F {\mathbb F}

\date{}

\title{On idempotents of a class of commutative rings}

\begin{document}

\author{ Fernanda D. de Melo Hern\'andez \thanks{Correspoding author}
\thanks{F. D. de Melo Hernandez (fdmelo@uem.br) is with the Departamento de Matem\'atica, Universidade Estadual de Maring\'a, Av. Colombo 5790, 87020-900, Maring\'a, PR, Brazil.} \and
C\'esar A. Hern\'andez Melo \thanks{C\'esar A. Hern\'andez Melo (cahmelo@uem.br) is with the Departamento de Matem\'atica, Universidade Estadual de Maring\'a, Av. Colombo 5790, 87020-900, Maring\'a, PR, Brazil.} 
\and Horacio Tapia-Recillas \thanks{Horacio Tapia-Recillas (htr@xanum.uam.mx) is with the Departamento de Matem\' aticas, Universidad Aut\' onoma Metropolitana-Unidad Iztapalapa, DF M\' exico.}}

\maketitle


\begin{abstract}
In the present work, a procedure for determining idempotents of a commutative ring having a sequence of ideals with certain properties is presented. As an application of this procedure, idempotent elements of various commutative rings are determined. Several examples are included illustrating the main results.


\end{abstract}


\section{Introduction} 

\quad Idempotent elements of an algebra are a topic of considerable research with a variety of applications including (theoretical) physics and chemistry (\cite{do},\cite{h},\cite{z},\linebreak \cite{l}), econometrics and various areas of mathematics, such as representation theory involving the decompositions of modules \cite{broche}. Information theory is not the exception, and idempotent elements are of great importance, particularly in coding theory with error detecting-correcting linear codes the alphabet of which is a finite (commutative) ring (\cite{bakshi},\cite{cruz},\cite{idemp},\cite{sergio},\cite{vera},\cite{fer},\cite{arora}).

Determining the idempotent elements of an arbitrary ring is not an easy task in general. In an attempt to give an answer to this problem, several works appear in the literature. For example
a set of primitive central idempotent elements are determined for the group algebra $\mathbb CG$ by using the group of characters of the finite group $G$ (\cite{livro}, Thm. 5.1.11, pag.185), this method can be interpreted in terms of the Discrete Fourier Transform.
In \cite{arora}, from the point of view of polynomial rings, the set of primitive central idempotent elements of $\mathbb F_q[x]/\langle x^n-1 \rangle$, where $\mathbb F_q$ is a finite field with $q$ elements, are determined for certain $q$ and $n$ in terms of cyclotomic classes. Considering a group algebra $\F_qG$ with G an abelian p-group in \cite{idemp}, under certain conditions on p and q, the set of primitive central idempotent elements is determined via the lattice of subgroups of the group $G$, extending the results of  \cite{arora}. In \cite{sergio} (see also \cite{vera}) the idempotents of the algebra $\Z_{p^k}[x]/\langle x^n-1\rangle$ are obtained by means of the irreducible factors of $x^n-1$ and the knowledge of idempotent elements in $\Z_{p}[x]/\langle x^n-1\rangle$, where $\Z_{p^k}$ is the ring of integers modulo $p^k$ with $p$ a prime and $k \geq 1$ an integer. 

Another way to determine idempotent elements in a ring $R$ is by ``lifting" idempotents elements. In \cite[Proposition 7.14, pag.405]{jacobson} idempotent elements of a ring $R$ are obtained by ``lifting" idempotent elements from the quotient ring $R/N$, where $N$ is a nil ideal of $R$.  
In this manuscript, starting from this result, 
under certain general conditions, a procedure to determine the idempotent elements of a commutative ring as ``lifted" idempotent elements of another ring
is given (Theorem \ref{IdemGeral}). 
This result has several consequences, for instance the set of idempotent elements are determined in cases which include commutative rings containing a nilpotent ideal; commutative group rings $RG$, where the ring $R$ contains a nilpotent ideal; commutative group rings $RG$ where the ring $R$ is a chain ring; and the commutative group ring $\Z_mG$.

The paper is organized as follows. In Section 2, basic facts about idempotent elements and, for completeness, the proof of {\cite[Proposition 7.14, pag.405]{jacobson}} is included. In Section 3, 
the main result of this manuscript, Theorem \ref{IdemGeral}, is presented. As a consequence, 
the set of idempotent elements of several rings is determined in Section 4, and
examples are included illustrating the main ideas. 


\section{Basic facts}

\quad The starting point for the results presented in this manuscript is the following appearing in {\cite[Proposition 7.14, pag.405]{jacobson}} on the construction of idempotent elements of a ring from those of a quotient ring. This result is recalled in the case of a commutative ring and for completeness, the main steps of the proof are included. The interested reader can see further details in the mentioned reference. 

Recall that an element $e$ of a ring $R$ is called {\it idempotent} if $e^2=e$, and two idempotent elements $e_1$ and $e_2$ of a ring $R$ are said to be {\it orthogonal} if $e_1 e_2 = e_2e_1 = 0$. An idempotent $e \in R$ is {\it primitive} if it can not be written as a sum of two non-trivial orthogonal idempotent elements.

\begin{propo}{\em \cite[Proposition 7.14]{jacobson}}\label{jacobson}
Let $R$ be a ring, $N$ a nil ideal of $R$ and $\bar{f} = f + N$ an idempotent element of the quotient ring $ R/N$. Then there exists an idempotent element $e$ in $R$ such that $\bar{e} = \bar{f}$, where `` $\bar{  }$ " denotes the canonical homomorphism from $R$ to $R/N$. Furthermore, if $R$ is commutative, the element $e$ is unique. 
\end{propo}

\proof
Since $\bar{f}$ is idempotent, $f^2-f \in N$ and because $N$ is a nil ideal, $(f^2-f)^n=0$ for some integer $n>0$. Then, if $g=1-f$, $0=(fg)^n=f^ng^n$. From the relation $f+g=1$ it follows that 
\[
1=1^{2n-1}=(f+g)^{2n-1}=h+e,
\]
where
\[
h=\sum_{i=0}^{n-1}\binom{2n-1}{i}f^ig^{2n-1-i}, \;\; e=\sum_{i=n}^{2n-1}\binom{2n-1}{i}f^ig^{2n-1-i}.
\]
Since $f^ng^n=0$ it follows that $eh=he=0$ and since $e+h=1$, $e^2=e$ and $h^2=h$. From this it is easy to see that $f^{2n-1} \equiv e $ mod $N$ and that $f \equiv f^2 \equiv \cdots \equiv f^{2n-1}$ mod $N$ from which it is concluded that $e \equiv f$ mod $N$. 

To prove uniqueness, consider the idempotent element of the form $e + z$ with $z$ nilpotent. The condition $(e+z)^2 =e+z$ gives $(1-2e)z= z^2$. Then $z^3 = (1-2e)z^2 = (1-2e)^2z$ and by induction, $(1-2e)^nz = z^{n+1}$. Since $(1-2e)^2 = 1-4e+4e = 1$, this implies that $z = 0$ and hence $e+z = e$. 
\cqd

Given an idempotent element $\overline {f}$ of $R/N$, if $R$ is commutative, the unique element $e \in R$ determined as in Proposition \ref{jacobson} will be called the {\it lifted} idempotent of $\overline {f}$.

\begin{obs}\label{Ob1}
The following observations are easy consequences of the previous result.
 \begin{enumerate}
 
 \item If $\bar f_1$ and $\bar f_2$ are orthogonal idempotent elements of $R/N$, then the corresponding lifted idempotent $e_1$ and $e_2$ of $R$ are also orthogonal. This follows from the fact that $e_1e_2$ is an idempotent element in the nil ideal $N$.
 
\item If $\bar f \in R/N$ is a primitive idempotent, then the corresponding lifted idempotent $e \in R$ is also primitive. In fact, if $e$ is not primitive, orthogonal idempotent elements $g, h$ exist in $R$ such that $e=g+h$. Then $\bar e=\bar f=\bar g + \bar h$ with $\bar g$ and $\bar h$ orthogonal idempotents in $R/N$, so either $\bar g=0$ or $\bar h=0$, i.e., $g \in N$ or $h \in N$ and the claim follows from the fact that $N$ is a nil ideal.

\item From the previous claims it follows that if $\{\bar f_1,...,\bar f_r \}$ is a set of primitive orthogonal idempotent elements in $R/N$, the corresponding set $\{e_1,...,e_r \}$ of lifted idempotent elements of $R$ has the same properties.
 
\item For a set $X$, $|X|$ will denote its cardinality. For a commutative ring $R$, let $E(R)$ denote the set of idempotent elements of $R$. Assuming the conditions of Proposition \ref{jacobson},
 \[
 |E(R)| = |E(R/N)|.
 \] 
 This follows from the fact that the canonical homomorphism from $R$ onto $R/N$ restricted to $E(R)$ is a bijection onto $E(R/N)$. The surjectivity follows from the construction of the idempotent element of $R$ from an idempotent element of $R/N$, and the injectivity follows from the uniqueness of this construction.

\end{enumerate}
\end{obs} 


\section{A method to compute lifted idempotents }

\quad In this section, under certain general conditions, a simple method to compute the lifted idempotent $e$ of $\bar{f}$ given in Proposition \ref{jacobson} is presented. More precisely, if a collection $\{N_1,\dots,N_k\}$ of ideals of a ring $R$ satisfies the CNC-condition (see Definition \ref{PosLiftIde} below), the lifted idempotent $e$ of the idempotent $f+N_1\in R/N_1$ can be computed as a power of $f$ in the ring $R$ (see Theorem \ref{IdemGeral} below). In order to do that the following result will be useful.

\begin{propo} \label{potencia}
Let $R$ be a commutative ring and $N$ a nilpotent ideal of index $t \geq 2$ in $R$. If $\bar f$ is an idempotent element in the quotient ring $R/N$ and $e$ is the corresponding lifted idempotent element in $R$ determined as in Proposition \ref{jacobson}, then:
\begin{enumerate}
	\item\label{pot} For any prime number $p$ such that $p \geq t$ and for all $n\in N$, there exists $r\in R$ such that $$(e+n)^p = e + pnr.$$ 
	\item\label{port} If a natural number $s>1$ exists, such that $sN=0$, and all the prime factors of the number $s$ are greater than or equal to the nilpotency index $t$ of the ideal $N$, then the lifted idempoted $e$ is
	$$e = f^s.$$ 
\item In particular, when the nilpotency index of the ideal $N$ is $t=2$ and $sN = 0$ for some $s\geq 2$, then $e = f^s$.
\end{enumerate}
\end{propo}

\proof \begin{enumerate}
\item  Since $n^t = 0$ and $e$ is an idempotent element in the ring $R$,
$$(e + n)^p = \sum_{j=0}^{p}\binom{p}{j}e^{p-j}n^j = e + \sum_{j=1}^{t-1}\binom{p}{j}en^j.$$
Since $p$ is a prime number, $p$ divides $\binom{p}{j}$ for all $1\leq j\leq p-1$.  Also, since $t\leq p$,
$$(e+n)^p=e+pn\left(k_1e+{k_2}en+\cdots + k_{t-1}en^{t-2}\right)$$
where $k_i=\tbinom{p}{i}/p$. Therefore,
$$(e+n)^{p}= e+pnr,$$
with $r=k_1e+k_2en+\cdots + k_{t-1}en^{t-2}\in R$.

\item Let $\{p_1,p_2,\ldots, p_m\}$ be the set of primes in the prime decomposition
of the integer $s$. 
Since $\bar f=\bar e$,  $f=e+n$ for some $n \in N$, and since $p_1\geq t$, from item \ref{pot} of Observation \ref{Ob1}, it follows that there exists $r_1\in R$ such that
$$f^{p_1} = (e + n)^{p_1}= e + p_1nr_1.$$
Similarly, since $p_2\geq t$ and $p_1nr_1\in N$, item \ref{pot} of Observation \ref{Ob1} implies that there exists $r_2\in R$ such that 
$$f^{p_1p_2}=(f^{p_1})^{p_2}=(e+p_1nr_1)^{p_2}=e+p_2(p_1nr_1)r_2.$$
Continuing with this process $r_3, r_4,\ldots, r_m\in R$ exists such that  
$$f^{s}=e+sn(r_1r_2\cdots r_m).$$ 
In other words, 
$$f^{s}=e+sh,$$ 
where $h=nr_1r_2\cdots r_s\in N$. Finally, since $h\in N$, and $sN=0$, it follows that  $e=f^{s}$. 

\item Observe that if the nilpotency index of the ideal $N$ is $t = 2$, it is obvious that all prime factors of $s$ are greater or equal to $t=2$. Therefore the proof of this claim is an immediate consequence of  item \ref{port} of Observation \ref{Ob1}.   
\end{enumerate}
\cqd

The following definition will play an important role in the rest of the manuscript.

\medskip
Recall that given an ideal $N$ of the ring $R$ and $k>1$, $N^k$ denotes the ideal generated by all products $x_1x_2 \cdots x_k$, where each $x_i \in N$ for $i= 1,2,...,k$.


\begin{defin}\label{PosLiftIde}
We say that a collection $\{N_1,..., N_k\}$ of ideals of a commutative ring $R$ satisfies the {\it CNC-condition} if the following properties hold:

\begin{enumerate}
\item \label{chc} {\bf Chain condition:} $\{0\}=N_k\subset N_{k-1}\subset\cdots \subset N_{2}\subset N_{1}\subset R$.

\item \label{nic} {\bf Nilpotency condition:} for $i=1,2,3,\ldots,k-1$, there exists an integer $t_i \geq 2$ such that $N_i^{t_i}\subset N_{i+1}$.

\item \label{cac} {\bf Characteristic condition:} for $i=1,2,3,\ldots,k-1$, there exists an integer $s_i\geq 1$ such that  $s_iN_i\subset N_{i+1}$. In addition, the prime factors of  $s_i$ are greater or equal to $t_i$.   
\end{enumerate}

The minimum number $t_i$ satisfying the nilpotency condition will be called the nilpotency index of the ideal $N_i$ in the ideal $N_{i+1}$. Similarly, the minimum number $s_i$ satisfying the characteristic condition will be called the characteristic of the ideal $N_i$ in the ideal $N_{i+1}$. 

\end{defin}

The nilpotency and characteristic conditions introduced above can be stated as follows:

\begin{itemize}
\item[a.]\label{pronic} The nilpotency condition is equivalent to the following: for $i=1,2,\ldots,k-1$, the quotient $N_{i}/N_{i+1}$ is a nilpotent ideal of index $t_i$ in the ring $R/N_{i+1}$.  In fact, if condition \ref{nic} is assumed to hold, then for all $n_1,n_2,\dots, n_{t_i} \in N_i$,  $n_1n_2n_3\cdots n_{t_i} \in N_{i}^{t_i}\subset N_{i+1}$. Thus,
$$
(n_1+N_{i+1})(n_2+N_{i+1})\cdots (n_{t_i}+N_{i+1})=n_1n_2\cdots n_{t_i}+N_{i+1}=N_{i+1},
$$ which implies that $N_{i}/N_{i+1}$ is a nilpotent ideal of index $t_i$ in the ring $R/N_{i+1}$. Conversely, if $N_{i}/N_{i+1}$ is a nilpotent ideal of index $t_i$ in the ring $R/N_{i+1}$, then for all $n_1,n_2,\dots, n_{t_i} \in N_i$, 
$$
(n_1+N_{i+1})(n_2+N_{i+1})\cdots (n_{t_i}+N_{i+1})=N_{i+1},
$$
which implies that $n_1n_2\cdots n_{t_i}\in N_{i+1}$. Thus any product of $t_i$ elements in the ideal $N_i$ are in the ideal $N_{i+1}$, hence $N_i^{t_i}\subset N_{i+1}$. Therefore, condition \ref{nic} holds.
   
\item[b.]\label{procac} The characteristic condition is equivalent to the following: for $i=1,2,\ldots, \linebreak k-1$, there exists a natural number $s_i\geq1$ such that $s_i(N_{i}/N_{i+1})=0$ in the ring $R/N_{i+1}$. In fact, if condition \ref{cac} holds, for all $n\in N_i$, $s_in\in N_{i+1}$. Then
$$
s_i(n+N_{i+1})=s_in+N_{i+1}=N_{i+1},
$$ 
implying that  $s_i(N_{i}/N_{i+1})=0$ in the ring $R/N_{i+1}$. Conversely, if we assume that $s_i(N_{i}/N_{i+1})=0$ in the ring $R/N_{i+1}$, then  for all $n\in N_i$, 
$$
s_i(n+N_{i+1})=N_{i+1},
$$ 
which implies that $s_in\in N_{i+1}$, proving that $s_iN_i\subset N_{i+1}$. Note that the characteristic $s_i$ of the ideal $N_i$ in the ideal $N_{i+1}$ satisfies $s_i\leq r_i$, where $r_i$ denotes the characteristic of the ring $R/N_{i+1}$. In fact, since for all $x\in R$, $r_ix+N_{i+1}=N_{i+1}$, $r_ix\in N_{i+1}$. Hence, $r_iR\subset N_{i+1}$, and then $r_iN_i\subset N_{i+1}$, showing that $s_i\leq r_i$.
\end{itemize}

Now, we are in a position to give the following result. 

\begin{theo}\label{IdemGeral} 
Let $	R$ be a commutative ring, $\{N_1, N_2, \ldots, N_k\}$ be a collection of ideals of $R$ satisfying the CNC-condition. If
$s_i$ is the characteristic of the ideal $N_i$ in the ideal $N_{i+1}$ and $f+N_1$ is an idempotent element of the ring $R/N_1$, then 
$$f^{s_1s_2\cdots s_{k-1}},$$
is an idempotent element of the ring $R$. Furthermore, $|E(R)|=|E(R/N_{k-1})|=\cdots=|E(R/N_{1})|.$
\end{theo} 

\proof  If $f+N_1$ is an idempotent element of the ring 
\begin{equation}\label{iso1}
R/N_1\cong \frac{(R/N_2)}{(N_1/N_2)},
\end{equation}
then $(f+N_2)+(N_1/N_2)$ is an idempotent element of the ring $\frac{(R/N_2)}{N_1/N_2}$. Since $N_{1}/N_{2}$ is a nilpotent ideal of index $t_1$ in the ring $R/N_{2}$ and $s_1$ satisfies the hypothesis of claim 2 of Proposition \ref{potencia}, it follows that  
$f^{s_1}+N_2$ is an idempotent element of the ring $R/N_2$. From the isomorphism
\begin{equation}\label{iso2}
R/N_2\cong \frac{(R/N_3)}{(N_2/N_3)},
\end{equation}
it follows that $(f^{s_1}+N_3)+(N_2/N_3)$ is an idempotent element of the ring $\frac{(R/N_3)}{N_2/N_3}$. Since $N_{2}/N_{3}$ is a nilpotent ideal of index $t_2$ in the ring $R/N_{3}$ and $s_2$ satisfies the hypothesis of claim 2  
of Proposition \ref{potencia}, $f^{s_1s_2}+N_3$ is an idempotent element of the ring $R/N_3$.   Continuing with this process, since 
\begin{equation}\label{isogen}
R/N_i\cong \frac{(R/N_{i+1})}{(N_i/N_{i+1})},
\end{equation}
$f^{s_1s_2\cdots s_i}+N_{i+1}$ is an idempotent element of the ring $R/N_{i+1}$. Finally, in the last step of the chain of ideals,  
$$
f^{s_1s_2\cdots s_{k-1}}+N_k=f^{s_1s_2\cdots s_{k-1}},
$$ 
is an idempotent element of the ring $R/N_k=R$.\\ 

From the ring isomorphism given in (\ref{isogen}), it follows that 
$$
|E(R/N_i)|=\left|E\left(\frac{R/N_{i+1}}{N_i/N_{i+1}}\right)\right|,
$$ 
for all $i=1,2,3\ldots,k-1.$ Since $N_i/N_{i+1}$ is a nilpotent ideal of index $t_i$ in $R/N_{i+1}$, and $R$ is a commutative ring by Remark \ref{Ob1},
$$
\left|E\left(\frac{R/N_{i+1}}{N_i/N_{i+1}}\right)\right|=|E(R/N_{i+1})|.
$$
Therefore $|E(R/N_{i})|=|E(R/N_{i+1})|,$ for all $i=1,2,3\ldots,k-1$, and the claim is proved.
\cqd

\begin{obs}
Observe that if $\{N_1,N_2,\ldots,N_k\}$ is a collection of ideals of the commutative ring $R$ satisfying the CNC-condition, any idempotent element $f+N_1$ of the ring $R/N_1$ is lifted up to the idempotent element $f^{s_1}+N_2$ of the ring $R/N_2$. This new idempotent element is lifted up to the idempotent element $f^{s_1s_2}+N_3$ of the ring $R/N_3$, and so on. At the end of this process, $f^{s_1s_2\cdots s_{k-1}}$ is an idempotent element of the ring $R$. The following chain of ring homomorphisms 
$$
R\xrightarrow{\phi_{k-1}}\frac{R}{N_{k-1}}\xrightarrow{\phi_{k-2}} \cdots \xrightarrow{\phi_{3}}\frac{R}{N_{3}}\xrightarrow{\phi_{2}}\frac{R}{N_{2}}\xrightarrow{\phi_{1}}\frac{R}{N_{1}},
$$
appears naturally in the lifting process of the idempotent element $f+N_1 \in R/N_1$. In addition, each homomorphism $\phi_i$ of this chain induces a bijection when restricted to the set of idempotent elements $E(R/N_{i+1})$, i.e., 
$$
{\phi_i}_{E}:E(R/N_{i+1})\rightarrow E(R/N_i),\hspace{0.7cm}x+N_{i+1}\rightarrow x +N_i  
$$
is bijective for $i=1,2,...,k-1$.
\end{obs}

\begin{obs}
Assuming the hypothesis of Theorem \ref{IdemGeral}, the set of idempotent elements $E(R)$ of the commutative ring $R$ is given by
$$E(R) = \{f^{s_1s_2\cdots s_{k-1}}:\,\,\, \bar{f} \in E(R/N_1) \}.$$
Hence, in order to determine the idempotent elements of the ring $R$, it is necessary to have a collection $\{N_1,N_2, \ldots, N_k\}$ of ideals of the ring $R$ satisfying the CNC-condition 
together with all idempotent elements of the quotient ring $R/N_1$.
\end{obs}


\section{Consequences and applications}

\quad In this section, by means of Theorem \ref{IdemGeral} the set of idempotent elements are determined in varios cases, including: commutative rings containing a nilpotent ideal; commutative group rings $RG$, where $R$ contains a nilpotent ideal; commutative group rings $RG$ where $R$ is a chain ring; the group ring $\Z_mG$, where $\Z_m$ is the ring of integers modulo $m$. In each of these cases examples are included illustrating the results.

\subsection{Commutative rings containing a nilpotent ideal}

The following result describes the set of idempotent elements of a commutative ring $R$ containing a nilpotent ideal $N$ in terms of the idempotent elements of the quotient ring $R/N$.

\begin{propo}\label{GeralNil} 
Let $R$ be a commutative ring and $N$ a nilpotent ideal of $R$ of nilpotency index $k\geq 2$.  Let $s>1$ be the characteristic of the quotient ring $R/N$. If $f+ N$ is an idempotent element of $R/N$, then 
$$f^{s^{k-1}}$$
is an idempotent element of the ring $R$. Moreover, $|E(R)|=|E(R/N)|$.
\end{propo} 
\proof The proof of this proposition is a consequence of Theorem \ref{IdemGeral}. It will be enough to show that the collection  $B=\{N, N^2,...,N^k\}$ of ideals of the ring $R$ satisfies the CNC-condition with  nilpotency index and characteristic of the ideal $N^{i}$ in the ideal $N^{i+1}$ being $t_i=2$ and $s_i=s$ for all $i=1,2,3,\dots,k-1$. Indeed, 
\begin{enumerate}
\item It is clear that the collection $B$ satisfies the chain condition. 
\item Since $(N^i)^2 = N^{2i}$ and $i+1\leq 2i$ for $i=1,2,3,\dots, k-1$, it follows that $(N^i)^2 \subset N^{i+1}$. Hence, the collection $B$ satisfies the nilpotency condition. 
\item Since the ring $R/N$ has characteristic $s$, there exists $n\in N$ such that $\sum_{i=1}^s1_R=n$. Then 
\begin{equation}
sN^{i} =(1_R+\cdots+1_R)N^{i}=nN^{i}\subset N^{i+1}, 
\end{equation}
and it follows that $sN^i\subset N^{i+1}$ for $i=1,2,3,\dots,k-1$. In addition, all prime factors of $s_i=s$ are greater or equal to the nilpotency index $t_i=2$, thus proving that the collection $B$ satisfies the characteristic condition.
\end{enumerate}
Therefore, Theorem \ref{IdemGeral} implies that $f^{s^{k-1}}$ is an idempotent element of the ring $R$ and  $|E(R)|=|E(R/N)|$. 
\cqd

An immediate consequence of Proposition \ref{GeralNil} is the following:

\begin{coro}\label{coroe}
Let $	R$ be a commutative ring, $a \in R$ a nilpotent element of index $k$, and $s>1$ the characteristic of the quotient ring $R/\langle a\rangle$. If $f+ \langle a \rangle$ is an idempotent element of $R/\langle a\rangle$ ,then 
$$f^{s^{k-1}}$$
is an idempotent element of the ring $R$. Moreover, $|E(R)|=|E(R/\langle a\rangle)|.$
\end{coro} 
\proof Since $N = \langle a \rangle$ is a nilpotent ideal of nilpotency index $k$ in $R$, the result follows from Proposition \ref{GeralNil}
\medskip
\cqd

The following example illustrates the previous corollary.

\begin{example}
Let 
$$
\mathbb{Z}_{p^{k}}[i]=\{a+bi: a,b\in \mathbb{Z}_{p^k}, i^2=-1\},
$$ 
where $p>2$ is a prime and $k>1$ a natural number. 
It is easy to see that $a=p$ is a nilpotent element of index $k$ in the ring $\mathbb{Z}_{p^{k}}[i]$. Since 
$$
\frac{\mathbb{Z}_{p^{k}}[i]}{\langle p\rangle}\cong \mathbb{Z}_{p}[i],
$$
and the ring $\mathbb{Z}_{p}[i]$ has characteristic $s=p$, Corollary \ref{coroe} implies that 
$$
E\left(\mathbb{Z}_{p^{k}}[i]\right)=\left\{f^{r}: r=p^{k-1},\text{ and }  \bar{f}\in E\left(\mathbb{Z}_{p}[i]\right)\right\}, 
$$
and
$$
\left|E\left(\mathbb{Z}_{p^{k}}[i]\right)\right|=\left|E\left(\mathbb{Z}_{p}[i]\right)\right|.
$$
The idempotent elements of the ring $\Z_p[i]$ are determined as follows. Since the idempotent elements of $\Z_p[i]$ are of the form $z = a + bi$, where $a=a^2-b^2$ and $b=2ab$, if $b = 0$ then $z = 0$ or $z = 1$. If $b \neq 0 $, $2a=1$ in $\Z_p$, hence $a=(p+1)/2$. From the relation $a=a^2-b^2$ it can be seen that $4b^2= (2b)^2=-1$ in $\Z_p$. Now we assume that $p \equiv 1 \mod 4$. Then $-1$ 
is a quadratic residue in $\mathbb Z_p$, (Theorem 3.1, pag.132, \cite{n-z}) and it is easy to see (by Wilson's Theorem) that the solutions are $2b=\pm (\frac{p-1}{2})!$. If $x=(\frac{p-1}{2})!$, the idempotent elements of the ring $\Z_p[i]$ are:

\begin{equation}\label{idemgauss}
\begin{aligned}
\bar{f_1}&= 0 ,\hspace{0.5cm} \bar{f_2}=1,\\
\bar{f_3}&=(p+1)/2 + (x/2)i,\hspace{0.5cm} \bar{f_4}=(p+1)/2 - (x/2)i.
\end{aligned}
\end{equation}
Hence, the set of idempotents elements of the ring $\mathbb{Z}_{p^{k}}[i]$ is  
$$
E\left(\mathbb{Z}_{p^{k}}[i]\right)=\left\{f_1^{r}, f_2^{r}, f_3^{r}, f_4^{r}\right\},
$$ where $r=p^{k-1}$.
\medskip
Observe that the ring $\Z_{p^k}[i]$ is isomorphic to the ring $\mathbb Z_{p^k}[x]/\langle x^2+1 \rangle$.
\end{example}


\subsection{Commutative group rings}

In the following the set of idempotent elements for the group ring $RG$, where $R$ is a commutative ring containing a collection of ideals satisfying the CNC-condition and $G$ a commutative group, is determined.
\begin{propo}\label{GeralGR} 
Let $	R$ be a commutative ring and $G$ 
a commutative group. Let $\{N_1, N_2, \ldots, N_k\}$ be a collection of ideals of $R$ satisfying the CNC-condition. 
Let  $s_i$ be the characteristic of the ideal $N_i$ in the ideal $N_{i+1}$. If $f + N_1G$ is an idempotent element of the group ring $(R/N_1)G$, then $$f^{s_1s_2\cdots s_{k-1}}$$
is an idempotent element of the group ring $RG$. Furthermore, $|E(RG)|=|E((R/N_{1})G)|$.
\end{propo} 
\proof The proof of this proposition is a consequence of Theorem \ref{IdemGeral}. First, we show that the collection $B=\{N_1G, N_2G, \ldots, N_kG\}$ of ideals of the group ring $RG$ satisfies the CNC-condition with nilpotency index and characteristic of the ideal $N_iG$ in the ideal $N,_{i+1}G$, exactly the same nilpotency index and characteristic of the ideal $N_i$ in the ideal $N_{i+1}$. Indeed, 

\begin{enumerate}
\item It is clear that the collection $B$ satisfies the chain condition. 
\item If $t_i$ denotes the nilpotency index of the ideal $N_i$ in the ideal $N_{i+1}$, then $N_{i}^{t_i}\subset N_{i+1}$. It is not difficult to see that  $(N_iG)^{t_i}=N_i^{t_i}G$. Thus,
$$
(N_iG)^{t_i}=N_i^{t_i}G\subset N_{i+1}G,
 $$
proving that the collection $B$ satisfies the nilpotency condition.
\item Since $s_i$ is the characteristic of the ideal $N_i$ in the ideal $N_{i+1}$, then $s_iN_{i}\subset N_{i+1}$. It is clear that  $s_i(N_iG)=(s_iN_i)G$. Hence,
$$
s_i(N_iG)=(s_iN_i)G\subset N_{i+1}G.
$$

Now, since the collection $\{N_1,N_2,\dots,N_k\}$ satisfies the CNC-condition, it is obvious that all prime factors of the characteristic $s_i$ are greater or equal to the nilpotency index $t_i$ for all $i=1,2,3,\dots,k-1$, thus proving that the collection $B$ satisfies the characteristic condition.
\end{enumerate}

From Theorem \ref{IdemGeral} and the isomorphism
$$
\frac{RG}{N_1G}\cong \left(\frac{R}{N_1}\right)G,
$$
it follows that $f^{s_1s_2\cdots s_{k-1}}$ is an idempotent element of the group ring $RG$, and $|E(RG)|=|E((R/N_{1})G)|$.
\cqd

\begin{coro}\label{GeralNilGR} 
Let $	R$ be a commutative ring, $N$ a nilpotent ideal of index $k$ in $R$, $G$ 
a commutative group, and $s$ the characteristic of the quotient ring $R/N$. If $f + NG$ is an idempotent element in the group ring $(R/N)G$, then 
$$f^{s^{k-1}}$$
is an idempotent element in the group ring $RG$. Furthermore, $|E(RG)|=|E((R/N)G)|$.
\end{coro}
\proof 
The proof of Proposition \ref{GeralNil} shows that the collection $\{N, N^2,...,N^{k}\}$ of ideals of the ring $R$ satisfies the CNC-condition with constant characteristic $s_i=s$ for $i=1,2,3,\cdots, k-1$. The result follows from Proposition \ref{GeralGR}.
\cqd 

\begin{example}
Let $\F_p$ be the finite field of order $p$, $R=\F_p [X_1,...,X_m]/\langle X_1^p,...,X_m^p \rangle$ and $N=\langle x_1,...,x_m \rangle$ the ideal generated by $x_j = X_j + \langle X_1^p,...,X_m^p \rangle$. It is not difficult to see that $N$ has nilpotency index equal to $m(p-1)+1$ and that $R/N \simeq \F_p$. Then for $G$ a finite abelian group, Corollary \ref{GeralNilGR} implies that
$$E(RG) = \{f^{p^{m(p-1)}} \hspace{0.1cm}:\hspace{0.1cm} {\bar f} \in E(\F_p G)\}.$$
\end{example}

\begin{obs}
It is easy to see that the collection of ideals 
\[
\{{ N}, {N}^2,...,{N}^{m(p-1)},  {N}^{m(p-1)+1}=(0)\},
\]
satisfies condition (CNC). This sequence of ideals is important in the study of Reed-Muller codes over the finite field $\F_p$, $p$ a prime, since for an integer $\nu$ with $0 \leq \nu \leq m(p-1)$,
\[
{N}^{m(p-1)-\nu} = RM_{\nu}(m,p),
\]  
where $RM_{\nu}(m,p)$ is the Reed-Muller code of order $\nu$ over the alphabet $\mathbb F_p$.
Details can be found in \cite{and}.
\end{obs}

\begin{coro}\label{ElNil}
Let $	R$ be a commutative ring, $a$ be a nilpotent element of index $k$ in $R$, $G$ 
a commutative group and $s$ the characteristic of the quotient ring $R/\langle a\rangle$. If $f + \langle a \rangle G$ is an idempotent element of the group ring $(R/\langle a\rangle)G$, then 
$$f^{s^{k-1}}$$
is an idempotent element of the group ring $RG$. Furthermore, $|E(RG)|=|E((R/\langle a\rangle)G)|$.
\end{coro} 
\proof 
It is enough to observe that $N = \langle a \rangle $ is a nilpotent ideal of index $k$ in $R$. The result follows from Corollary \ref{GeralNilGR}.
\cqd

\subsection{Commutative group ring $RG$ with $R$ a chain ring}

Let $R$ be a finite commutative chain ring and $G$ 
a commutative group. It is well known that $R$ contains a unique maximal nilpotent ideal $N=\langle a\rangle$ for some $a\in R$. If $k$ denotes the nilpotency index of $a$, and $p$ the characteristic of the (residue) field $\F=R/\langle a\rangle$ from Corollary \ref{ElNil}, it follows that  
\begin{equation}\label{few}
E(RG)=\left\{f^{r}: r=p^{k-1},\hspace{0.1cm}\text{and}\hspace{0.1cm} \bar{f}\in E(\F G)\right\}.
\end{equation}
Examples of finite commutative chain rings include the modular ring of integers $R=\mathbb{Z}_{p^k}$, where $p$ is a prime number and $k>1$ is an integer. In this example, the maximal nilpotent ideal is $N=\langle p \rangle$, $p$ has nilpotency index equal to $k$ in $\mathbb{Z}_{p^k}$, and $\F\cong \Z_p$. Thus, 
\begin{equation}\label{few1}
E(\mathbb{Z}_{p^k}G)=\left\{f^{r}: r=p^{k-1},\hspace{0.1cm}\text{and}\hspace{0.1cm} \bar{f}\in E(\mathbb{Z}_pG)\right\}.
\end{equation}
When the group $G$ is cyclic of order $n$, it is known that $\mathbb{Z}_{p^k}G\cong \mathbb{Z}_{p^k}[x]/\langle x^n-1\rangle$, then (\ref{few1}) can be written in the following equivalent form:
$$
E(\Z_{p^k}[x]/\langle x^n-1 \rangle)=\left\{f^{r}: r=p^{k-1},\hspace{0.1cm}\text{and}\hspace{0.1cm} \bar{f}\in E(\Z_p[x]/\langle x^n-1 \rangle)\right\}.
$$
 
Another interesting class of chain rings are Galois rings, $R=\mathbb{Z}_{p^k}[x]/\langle q(x) \rangle$ where $p$ is a prime number, $k>1$ and $q(x)$ is a monic polynomial of degree $r$ image of which in $\Z_{p}[x]$ is irreducible. In this example, the maximal nilpotent ideal is $N=pR$, $p$ has nilpotency index equal to $k$ in $R$, and $\F=R/N\cong \F_{p^{r}}$ is the finite field of order $p^{r}$. Thus, 
\begin{equation}\label{few2}
E((\mathbb{Z}_{p^k}[x]/\langle q(x) \rangle)G)=\left\{f^{r}: r=p^{k-1},\hspace{0.1cm}\text{and}\hspace{0.1cm} \bar{f}\in E(\F_{p^{r}}G)\right\}.
\end{equation}

Combining the theory of Ferraz and Polcino in \cite{idemp} with the results discussed above, it will be shown how to compute the idempotent elements in group rings $\Z_{p^k}G$ for certain commutative groups $G$. This is illustrated in the following examples.

\begin{example}
Let $p=5$ and let $G=C_7=\langle g \rangle =\{e, g, g^2,g^3,g^4,g^5,g^6\}$ be the cyclic group of order 7 generated by $g$. From {\cite[Lemma 3]{idemp}} and {\cite[Corollary 4]{idemp}} it can be seen that the idempotent elements of the group algebra $\mathbb{Z}_5C_7$ are: 
\begin{equation*}
\begin{aligned}
\bar{f_1}&= 0 ,\hspace{0.5cm} \bar{f_2}=1,\\
\bar{f_3}&= \widehat{G} = 3e+3g+3g^2 +3g^3+3g^4+3g^5+3g^6, \\
\bar{f_4}&= 1 - \widehat{G} = 3e+2g+2g^2 +2g^3+2g^4+2g^5+2g^6,
 \end{aligned}
\end{equation*}
where for a subgroup $H$ of $G$, the idempotent element $\dis{\widehat{H}}$ is defined by 
$$
\dis{\widehat{H}=\frac{1}{|H|} \sum_{h \in H} h}.
$$
Therefore, from (\ref{few}), the idempotent elements of the group ring $\mathbb{Z}_{5^3}C_7$ are:
\begin{equation*}
\begin{aligned}
f_1^{5^2}&= 0 ,\hspace{0.5cm} f_2^{5^2}=1,\\
f_3^{5^2}&=18e+18g+18g^2 +18g^3+18g^4+18g^5+18g^6,\\
f_4^{5^2}&=108e+107g+107g^2 +107g^3+107g^4+107g^5+107g^6.
\end{aligned}
\end{equation*}
\end{example}

\begin{example} \label{exabeliano}
Let $p=2$ and $G=\langle a, b | a^5 = b^5 = 1, ab = ba \rangle$ be an abelian group of order $25$. From {\cite[Lemma 5]{idemp}} and {\cite[Theorem 4.1]{idemp}}, ({see also \cite[Section III]{fer}}), it follows that the primitive orthogonal idempotent elements of the group ring $\mathbb{Z}_2(C_5 \times C_5)$ are:
\begin{equation*}
\begin{array}{ll}
f_0 = \widehat{G}, & f_1= \widehat{\langle a \rangle} - \widehat{G},\\
f_2 = \widehat{\langle b \rangle} - \widehat{G},  & f_i=\widehat{\langle ab^i \rangle} - \widehat{G}, \,\, i =1,2,3,4,\\
\end{array}
\end{equation*}
where for a subgroup $H$ of $G$, $\dis{\widehat{H} = \frac{1}{|H|} \sum_{h \in H} h}.$

Therefore, the set of primitive orthogonal idempotent elements of the group ring $\mathbb{Z}_{2^3}(C_5 \times C_5)$ is
$$ E(\mathbb{Z}_{2^3}(C_5 \times C_5)) = \{f_0^{2^2}, f_1^{2^2}, f_2^{2^2}, f_i^{2^2} \},$$
and the other idempotents are obtained as a sum of the primitive orthogonal idempotent elements.

In general, for any $k \in \N$, the set of primitive orthogonal idempotent elements of the group ring $\mathbb{Z}_{2^k}(C_5 \times C_5)$ is
$$ E(\mathbb{Z}_{2^k}(C_5 \times C_5)) = \{f_0^{2^{k-1}}, f_1^{2^{k-1}}, f_2^{2^{k-1}}, f_i^{2^{k-1}} \}.$$
\end{example}


\subsection{Commutative group rings $\mathbb{Z}_mG$ with $G$ a commutative group}

The previously discussed results for $\Z_{p^k}G$ can be extended to the group ring $\Z_m G$, where $\Z_m$ is the ring of integers modulo $m>1$, and $G$ a finite commutative group.

\begin{theo}\label{Zn}
Let $m=p_1^{r_1}p_2^{r_2}\cdots p_j^{r_j}$ be the prime factorization of the integer $m \geq 2$. 
Set $m_i=m/p_i^{r_i}$ and let $s_i$ be a natural number such that $s_im_i=1\mod(p_i^{r_i})$ for $i=1,2,3,\cdots,j$. If $\bar{f_i}$ is an idempotent element of the group ring $\mathbb{Z}_{p_i}G$ for $i=1,2,3,\cdots,j$, by setting $\alpha_i=p_i^{r_i-1}$, the element $e$ given by
\begin{equation}\label{ide}
e=s_1m_1f_1^{\alpha_1}+s_2m_2f_2^{\alpha_2}+\cdots+s_jm_jf_j^{\alpha_j},
\end{equation}
is an idempotent element of the group ring $\mathbb{Z}_{m}G$. Moreover,
\begin{equation}\label{equip1}
|E(\mathbb{Z}_{m}G)|=|E(\mathbb{Z}_{p_1}G)||E(\mathbb{Z}_{p_2}G)|\cdots|E(\mathbb{Z}_{p_j}G)|.
\end{equation}
\end{theo}
\proof From the Chinese Remainder Theorem (CRT), $\Z_m\cong \Z_{p_1^{r_1}}\times \Z_{p_2^{r_2}}\times \cdots\times\Z_{p_j^{r_j}}$ and, therefore,
\begin{equation}\label{is}
\mathbb{Z}_{m}G\cong_{\phi}\mathbb{Z}_{p_1^{r_1}}G \times \mathbb{Z}_{p_2^{r_2}}G \times\cdots\times \mathbb{Z}_{p_j^{r_j}}G.
\end{equation}
If $\bar{f_i}$ is an idempotent element of $\mathbb{Z}_{p_i}G$, for $i=1,2,3,\cdots,j$, from Corollary \ref{ElNil},  $f_i^{\alpha_i}$ is an idempotent element of $\mathbb{Z}_{p_i^{r_i}}G$ for $i=1,2,3,\cdots,j$. Thus,  
$$
(f_1^{\alpha_1},f_2^{\alpha_2},\cdots,f_j^{\alpha_j}),
$$
is an idempotent element of the product ring given in (\ref{is}). Consequently, $\phi^{-1}\left( 
(f_1^{\alpha_1},f_2^{\alpha_2},\cdots,f_j^{\alpha_j})\right)=e$ is an idempotent element of the group ring $\mathbb{Z}_{m}G$. Finally, from the CRT, $e$ can be expressed in the form (\ref{ide}). The equality (\ref{equip1}) follows from Corollary \ref{ElNil} and (\ref{is}). 
\cqd

The following result provides an alternative way to compute the idempotents of $\Z_mG$.

\begin{theo}\label{Znv2}
Let $m=p_1^{r_1}p_2^{r_2}\cdots p_j^{r_j}$ be the prime factorization of the integer $m \geq 2$. 
Set $k=\max\{r_1,r_2,\ldots,r_j\}$, $c_i=(p_1p_2\cdots p_j)/p_i$, and let $t_i$ be a natural number such that $t_ic_i=1\mod(p_i)$ for $i=1,2,3,\cdots,j$. If $\bar{f_i}$ is an idempotent element of the group ring $\mathbb{Z}_{p_i}G$ for $i=1,2,3,\cdots,j$, then
\begin{equation}\label{idewv2}
e=[t_1c_1f_1+t_2c_2f_2+\cdots+t_jc_jf_j]^{(p_1p_2\cdots p_j)^{k-1}}
\end{equation}
is an idempotent element of the group ring $\mathbb{Z}_{m}G$. Therefore,
\begin{equation}\label{equipv2}
|E(\mathbb{Z}_{m}G)|=|E(\mathbb{Z}_{p_1p_2\cdots p_j}G)|.
\end{equation}
\end{theo}

\proof If $\bar{f_i}$ is an idempotent element of the group ring $\mathbb{Z}_{p_i}G$ for $i=1,2,3,\cdots,j$, from Theorem \ref{Zn} it follows that 
\begin{equation}\label{idfg}
g=t_1c_1\bar{f}_1+t_2c_2\bar{f}_2+\cdots+t_jc_j\bar{f}_j,
\end{equation}
is an idempotent element of the group ring $\mathbb{Z}_{p_1p_2\cdots p_j}G$. Observe that $a=p_1p_2\cdots p_j$ has nilpotency index $k$ in the ring $\mathbb{Z}_m$ and
$$
\frac{\mathbb{Z}_m}{a\mathbb{Z}_m}\cong \mathbb{Z}_a
$$
has characteristic $a=p_1p_2\cdots p_j$. Hence, from Corollary \ref{ElNil}, it follows that,
\begin{equation*}\label{idev2}
e=g^{(p_1p_2\cdots p_j)^{k-1}}
\end{equation*}
is an idempotent element of the group ring $\mathbb{Z}_{m}G$, proving (\ref{idewv2}). Relation (\ref{equipv2}) also follows from Corollary \ref{ElNil}.
\cqd

It is worth noting that Theorem \ref{Zn} can be extended to the group ring $RG$, where $G$ is a finite commutative group and $R = \F_q[x]/\langle m(x) \rangle$, with $\F_q$ a finite field with $q = p^s$ elements, and $m(x)\in\F_q[x]$. Indeed, let $m(x)=p_1^{r_1}(x)p_2^{r_2} (x)\cdots p_j^{r_j}(x)$ be the factorization of the polynomial $m(x)$ in $\F_q[x]$,  $m_i(x) = m(x)/p_i^{r_i}(x)$ and  $s_i(x)\in\F_q[x]$ such that $s_i(x)m_i(x) \equiv 1 \mod(p_i^{r_i}(x))$ for $i=1,2,\dots,j$. If $\bar{f_i}$ is an idempotent element in the group ring $\F_q[x]/\langle p_i(x) \rangle)G$, the element $e$ given by
\begin{equation}\label{ide1}
e=s_1(x)m_1(x)f_1^{\alpha_1}(x)+s_2(x)m_2(x)f_2^{\alpha_2}(x)+\cdots+s_j(x)m_j(x)f_j^{\alpha_j}(x),\hspace{0.3cm}\alpha_i=p^{r_i-1},
\end{equation}
is an idempotent element of the group ring $(\F_q[x]/\langle m(x) \rangle)G$. Furthermore,
\begin{equation}\label{equip}
|E(\F_q[x]/\langle m(x) \rangle G)|=|E((\F_q[x]/\langle p_1(x) \rangle)G)|\cdots|E((\F_q[x]/\langle p_j(x) \rangle)G)|.
\end{equation}

The proof of the former result is similar to the one given for Theorem \ref{Zn}. However, in the following lines the main steps of the proof are given. If $R_i=\F_q[x]/\langle p_i^{r_i}(x) \rangle$, from the Chinese Remainder Theorem (CRT), it follows that
$$
RG \cong_{\phi_1} R_1G \times \cdots\times R_j G.
$$
Since $a_i=p_i(x)+ \langle p_i^{r_i}(x)\rangle$ is a nilpotent element in the ring $R_i$ of nilpotency index $r_i$ and the field $R_i/\langle a_i \rangle \cong \F_q[x]/\langle p_i(x) \rangle$ has characteristic $p$ for all $i=1,2\dots,j$, from Corollary \ref{ElNil} it follows that, if $\bar{f_i}$ is an idempotent element in $(R_i/\langle a_i \rangle) G$, for each $i=1,2,3,\dots,j$, then $f_i^{\alpha_i}$ is an idempotent element in $R_iG$ for each $i=1,2,3,\dots,j$. Thus, $
(f_1^{\alpha_1},f_2^{\alpha_2},\cdots,f_j^{\alpha_j})$
is an idempotent element in the ring  $R_1G \times \cdots\times R_j G$. Consequently, $\phi_1^{-1}((f_1^{\alpha_1},f_2^{\alpha_2},\cdots,f_j^{\alpha_j}))=e$ is an idempotent element in the ring $RG$. Finally, from the CRT, $e$ can be expressed in the form (\ref{ide1}).


\medskip
Next, Theorems \ref{Zn} and \ref{Znv2} will be illustrated with examples. In example \ref{srt}, the idempotent elements of the group ring $\mathbb Z_{200}C_3\cong \mathbb Z_{200}[x]/\langle x^3-1\rangle$ are determined. In example \ref{srt2}, the idempotent elements of the group ring $\mathbb Z_{936}(C_5\times C_5)$ are given. Here, $C_n$ denotes the cyclic group of order $n$.
\begin{example}\label{srt}
Let $R=\mathbb{Z}_{200}C_3$, where $C_3=\langle g \rangle=\{g_0,g_1,g_2\}$ is the cyclic group of order three. First note that the  elements
\begin{equation}\label{ide2}
\begin{aligned}
u_1&=g_0+g_1+g_2,\hspace{2.5cm}u_3&=g_0,\\
u_2&=g_1+g_2,\hspace{3.4cm}u_4&=0,
\end{aligned}
\end{equation}
are the 4 idempotent elements of the group ring $\mathbb{Z}_{2}C_3$ (the calculation of these idempotent elements can be done by hand). In addition, 
 \begin{equation}\label{ide5}
\begin{aligned}
v_1&=4g_0+3g_1+3g_2,\hspace{2.1cm}v_3&=g_0,\\
v_2&=2g_0+2g_1+2g_2,\hspace{2.1cm}v_4&=0,
\end{aligned}
\end{equation}
are the 4 idempotent elements of the group ring  $\mathbb{Z}_{5}C_3$. In order to determine the idempotent elements of the group ring $R$, Theorem \ref{Zn} or  Theorem \ref{Znv2} are applied.
\begin{itemize}
\item[$\bullet$] In terms of Theorem \ref{Zn}: since $m=200=2^35^2$, $p_1=2$, $p_2=5$, $m_1=5^2$ and $m_2=2^3$, so $s_1=1$ and $s_2=22$. Then, 
\begin{equation}\label{ide200}
E(\mathbb{Z}_{200}C_3)=\{25f^{4}+176g^5: f\in E(\mathbb{Z}_{2}C_3),\hspace{0.2cm} g\in E(\mathbb{Z}_{5}C_3)\}.
\end{equation}
\item[$\bullet$] In terms of Theorem \ref{Znv2}: since $m=200=2^35^2$, $p_1=2$, $p_2=5$, $k=3$, $c_1=5$, and $c_2=2$ , so $t_1=1$ and $t_2=3$. Then, 
\begin{equation}\label{ide2002}
E(\mathbb{Z}_{200}C_3)=\{(5f+6g)^{100}: f\in E(\mathbb{Z}_{2}C_3),\hspace{0.2cm} g\in E(\mathbb{Z}_{5}C_3)\}.
\end{equation}
\end{itemize} 
From relations (\ref{ide2}), (\ref{ide5}) and by using an algorithm implemented by the authors in the programing language C, one can see that 
\begin{equation*}
\begin{aligned}
h_1&=184g_0+8g_1+8g_2,\hspace{2.5cm}h_{9}=126g_0+125g_1+125g_2,\\
h_2&=192g_0+192g_1+192g_2,\hspace{1.7cm}h_{10}=17g_0+192g_1+192g_2,\\
h_3&=142g_0+117g_1+117g_2,\hspace{1.7cm}h_{11}=67g_0+67g_1+67g_2,\\
h_4&=134g_0+133g_1+133g_2,\hspace{1.7cm}h_{12}=59g_0+83g_1+83g_2,\\
h_5&=51g_0+75g_1+75g_2,\hspace{2.3cm}h_{13}=9g_0+8g_1+8g_2,\\
h_6&=150g_1+125g_1+125g_2,\hspace{1.7cm}h_{14}=75g_0+75g_1+75g_2,\\
h_7&=g_0,\hspace{5.3cm}h_{15}=25g_0,\\
h_8&=0,\hspace{5.4cm}h_{16}=176g_0.
\end{aligned}
\end{equation*}
are the 16 idempotent elements of $\mathbb{Z}_{200}C_3$.
\end{example}

\begin{example}\label{srt2}
Let $Z_{936}G$ where $G= C_5 \times C_5=\langle a, b | a^5 = b^5 = 1, ab = ba \rangle $ is an abelian group of order $25$. First, note that from example \ref{exabeliano},
\begin{equation*}
\begin{array}{ll}
f_0 = \widehat{G}, & f_1= \widehat{\langle a \rangle} - \widehat{G},\\
f_2 = \widehat{\langle b \rangle} - \widehat{G},  & f_i=\widehat{\langle ab^i \rangle} - \widehat{G}, \,\, i =1,2,3,4,\\
\end{array}
\end{equation*}
are the primitive orthogonal idempotent elements of the group ring $\mathbb{Z}_2(C_5 \times C_5)$.

Also, from {\cite[Lemma 5]{idemp}}, {\cite[Theorem 4.1]{idemp}} and {\cite[Section III]{fer}} it follows that the primitive orthogonal idempotent elements of the group ring $\mathbb{Z}_3(C_5 \times C_5)$ are:
\begin{equation*}
\begin{array}{ll}
g_0 = \widehat{G}, & g_1= \widehat{\langle a \rangle} - \widehat{G},\\
g_2 = \widehat{\langle b \rangle} - \widehat{G},  & g_i=\widehat{\langle ab^i \rangle} - \widehat{G}, \,\, i =1,2,3,4,\\
\end{array}
\end{equation*}
and the primitive orthogonal idempotent elements of the group ring $\mathbb{Z}_{13}(C_5 \times C_5)$ are:
\begin{equation*}
\begin{array}{ll}
h_0 = \widehat{G}, & h_1= \widehat{\langle a \rangle} - \widehat{G},\\
h_2 = \widehat{\langle b \rangle} - \widehat{G},  & h_i=\widehat{\langle ab^i \rangle} - \widehat{G}, \,\, i =1,2,3,4.\\
\end{array}
\end{equation*}

In order to determine the idempotent elements of the group ring $Z_{936}G$, Theorem \ref{Zn} can be applied. Since $m=936=2^33^213$, $p_1=2$, $p_2=3$, $p_3 = 13$, $m_1=3^213$, $m_2=2^313$ and $m_3 = 2^33^2$, so $s_1=5$, $s_2 = 2$and $s_3 = 2$. Then, 
\begin{equation}\label{ide200}
E(\Z_{936}G)=\{585f^{4}+208g^3 + 144h: f\in E(\Z_{2}G),\hspace{0.2cm} g\in E(\Z_{3}G)\hspace{0.2cm} h\in E(\Z_{13}G)\}.
\end{equation}
is a set of primitive orthogonal idempotent elements of the group ring $Z_{936}G$.
\end{example}




\section*{Acknowledgements}
This research was carried on while the first author was a post-doctoral scholar at the Departamento de Matem\'aticas, Universidad Aut\'onoma Metropolitana - Iztapalapa, Cd. de M\'exico, M\'exico. The second author was a visiting researcher at the Instituto de Investigaciones en Matem\'aticas Aplicadas y Sistemas (IIMAS), Universidad Nacional Aut\'onoma de M\'exico, Cd. de M\'exico. They would like to thank these institutions for the warm reception during their stay.



\end{document}